\documentclass{article}
\usepackage{amssymb,amsmath,
theorem,euscript}

\input{epsf}

\usepackage[dvips]{graphicx}

\newcounter{sec}

\newcounter{punct}[sec]

\def\punct{\refstepcounter{punct}{\arabic{sec}.\arabic{punct}.  }}

\def\COUNTERS{\addtocounter{sec}{1}
              \setcounter{punct}{0}
          \setcounter{equation}{0}
          \setcounter{theorem}{0}
          \setcounter{problem}{0}
            \setcounter{Apunct}{0}
          }

\newtheorem{theorem}{Theorem}[sec]
\newtheorem{proposition}[theorem]{Proposition}

\def\COUNTERS{\addtocounter{sec}{1}
              \setcounter{punct}{0}
          \setcounter{equation}{0}
          \setcounter{theorem}{0}
          }

\begin{document}

\def\SL{\mathrm {SL}}
\def\SU{\mathrm {SU}}
\def\GL{\mathrm  {GL}}
\def\U{\mathrm  U}
\def\OO{\mathrm  O}
\def\Sp{\mathrm  {Sp}}
\def\SO{\mathrm  {SO}}
\def\SOS{\mathrm {SO}^*}

\def\PGL{\mathrm  {PGL}}
\def\PU{\mathrm {PU}}

\def\Gr{\mathrm{Gr}}

\def\Fl{\mathrm{Fl}}

\def\OSp{\mathrm {OSp}}

\def\Mat{\mathrm{Mat}}

\def\Pfaff{\mathrm {Pfaff}}

\def\Ind{\mathrm{Ind}}

\def\B{\mathbf B}

\def\phi{\varphi}
\def\epsilon{\varepsilon}
\def\kappa{\varkappa}

\def\le{\leqslant}
\def\ge{\geqslant}

\renewcommand{\Re}{\mathop{\rm Re}\nolimits}

\renewcommand{\Im}{\mathop{\rm Im}\nolimits}

\def\pia{\pi_\downarrow}

\newcommand{\im}{\mathop{\rm im}\nolimits}
\newcommand{\indef}{\mathop{\rm indef}\nolimits}
\newcommand{\dom}{\mathop{\rm dom}\nolimits}
\newcommand{\codim}{\mathop{\rm codim}\nolimits}

\def\cA{\mathcal A}
\def\cB{\mathcal B}
\def\cC{\mathcal C}
\def\cD{\mathcal D}
\def\cE{\mathcal E}
\def\cF{\mathcal F}
\def\cG{\mathcal G}
\def\cH{\mathcal H}
\def\cJ{\mathcal J}
\def\cI{\mathcal I}
\def\cK{\mathcal K}
\def\cL{\mathcal L}
\def\cM{\mathcal M}
\def\cN{\mathcal N}
\def\cO{\mathcal O}
\def\cP{\mathcal P}
\def\cQ{\mathcal Q}
\def\cR{\mathcal R}
\def\cS{\mathcal S}
\def\cT{\mathcal T}
\def\cU{\mathcal U}
\def\cV{\mathcal V}
\def\cW{\mathcal W}
\def\cX{\mathcal X}
\def\cY{\mathcal Y}
\def\cZ{\mathcal Z}

\def\frA{\mathfrak A}
\def\frB{\mathfrak B}
\def\frC{\mathfrak C}
\def\frD{\mathfrak D}
\def\frE{\mathfrak E}
\def\frF{\mathfrak F}
\def\frG{\mathfrak G}
\def\frH{\mathfrak H}
\def\frJ{\mathfrak J}
\def\frK{\mathfrak K}
\def\frL{\mathfrak L}
\def\frM{\mathfrak M}
\def\frN{\mathfrak N}
\def\frO{\mathfrak O}
\def\frP{\mathfrak P}
\def\frQ{\mathfrak Q}
\def\frR{\mathfrak R}
\def\frS{\mathfrak S}
\def\frT{\mathfrak T}
\def\frU{\mathfrak U}
\def\frV{\mathfrak V}
\def\frW{\mathfrak W}
\def\frX{\mathfrak X}
\def\frY{\mathfrak Y}
\def\frZ{\mathfrak Z}

\def\fra{\mathfrak a}
\def\frb{\mathfrak b}
\def\frc{\mathfrak c}
\def\frd{\mathfrak d}
\def\fre{\mathfrak e}
\def\frf{\mathfrak f}
\def\frg{\mathfrak g}
\def\frh{\mathfrak h}
\def\fri{\mathfrak i}
\def\frj{\mathfrak j}
\def\frk{\mathfrak k}
\def\frl{\mathfrak l}
\def\frm{\mathfrak m}
\def\frn{\mathfrak n}
\def\fro{\mathfrak o}
\def\frp{\mathfrak p}
\def\frq{\mathfrak q}
\def\frr{\mathfrak r}
\def\frs{\mathfrak s}
\def\frt{\mathfrak t}
\def\fru{\mathfrak u}
\def\frv{\mathfrak v}
\def\frw{\mathfrak w}
\def\frx{\mathfrak x}
\def\fry{\mathfrak y}
\def\frz{\mathfrak z}

\def\fros{\mathfrak{s}}

\def\bfa{\mathbf a}
\def\bfb{\mathbf b}
\def\bfc{\mathbf c}
\def\bfd{\mathbf d}
\def\bfe{\mathbf e}
\def\bff{\mathbf f}
\def\bfg{\mathbf g}
\def\bfh{\mathbf h}
\def\bfi{\mathbf i}
\def\bfj{\mathbf j}
\def\bfk{\mathbf k}
\def\bfl{\mathbf l}
\def\bfm{\mathbf m}
\def\bfn{\mathbf n}
\def\bfo{\mathbf o}
\def\bfp{\mathbf q}
\def\bfr{\mathbf r}
\def\bfs{\mathbf s}
\def\bft{\mathbf t}
\def\bfu{\mathbf u}
\def\bfv{\mathbf v}
\def\bfw{\mathbf w}
\def\bfx{\mathbf x}
\def\bfy{\mathbf y}
\def\bfz{\mathbf z}

\def\bfA{\mathbf A}
\def\bfB{\mathbf B}
\def\bfC{\mathbf C}
\def\bfD{\mathbf D}
\def\bfE{\mathbf E}
\def\bfF{\mathbf F}
\def\bfG{\mathbf G}
\def\bfH{\mathbf H}
\def\bfI{\mathbf I}
\def\bfJ{\mathbf J}
\def\bfK{\mathbf K}
\def\bfL{\mathbf L}
\def\bfM{\mathbf M}
\def\bfN{\mathbf N}
\def\bfO{\mathbf O}
\def\bfP{\mathbf P}
\def\bfQ{\mathbf Q}
\def\bfR{\mathbf R}
\def\bfS{\mathbf S}
\def\bfT{\mathbf T}
\def\bfU{\mathbf U}
\def\bfV{\mathbf V}
\def\bfW{\mathbf W}
\def\bfX{\mathbf X}
\def\bfY{\mathbf Y}
\def\bfZ{\mathbf Z}

\def\R {{\mathbb R }}
 \def\C {{\mathbb C }}
  \def\Z{{\mathbb Z}}
  \def\H{{\mathbb H}}
\def\K{{\mathbb K}}
\def\N{{\mathbb N}}
\def\Q{{\mathbb Q}}
\def\A{{\mathbb A}}

\def\T{\mathbb T}

\def\bbA{\mathbb A}
\def\bbB{\mathbb B}
\def\bbD{\mathbb D}
\def\bbE{\mathbb E}
\def\bbF{\mathbb F}
\def\bbG{\mathbb G}
\def\bbI{\mathbb I}
\def\bbJ{\mathbb J}
\def\bbL{\mathbb L}
\def\bbM{\mathbb M}
\def\bbN{\mathbb N}
\def\bbO{\mathbb O}
\def\bbP{\mathbb P}
\def\bbQ{\mathbb Q}
\def\bbS{\mathbb S}
\def\bbT{\mathbb T}
\def\bbU{\mathbb U}
\def\bbV{\mathbb V}
\def\bbW{\mathbb W}
\def\bbX{\mathbb X}
\def\bbY{\mathbb Y}

 \def\ov{\overline}
\def\wt{\widetilde}
\def\wh{\widehat}

\def\P{\mathbb P}

\def\bO{\bf O}

\def\arr{\rightrightarrows}

\def\SS{\smallskip}

\def\ev{{\mathrm{even}}}
\def\od{{\mathrm{odd}}}

\def\q{\quad}

\def\F{\mathbf F}

\def\b{\mathbf b}

\def\RA{\Longrightarrow}

\def\la{\langle}
\def\ra{\rangle}

\begin{center}
\bf\Large On spherical functions on the group
$\SU(2)\times\SU(2)\times\SU(2)$

\sc\large

\bigskip

Yury Neretin%
\footnote{Supported by the grant  FWF,project P19064,
the Russian Federal Agency on Nuclear Energy,
and grants,
 NWO.047.017.015, JSPS-RFBR-07.01.91209.}

\end{center}

{\small We consider the group 
$G:=\SU(2)\times\dots\times\SU(2)$
($l$ times), where $\SU(2)$ is the group of unitary
$2\times 2$ matrices, whose determinant is $1$.
Denote by   $K\simeq\SU(2)$ the diagonal subgroup in
 $G$.
We obtain a generating function for all 
 $K$-spherical functions on  $G$.}

\begin{flushright}
To G.~I.~Olshanski in his 60th birthday
\end{flushright}

\section{Introduction}

\COUNTERS

{\bf \punct The group $\SU(2)$.} Denote by 
 $\SU(2)$ the group of unitary $2\times 2$ matrices with determinant
$=1$; recall that any such matrix has the form
$$
g=\begin{pmatrix}a&b\\ -\ov b&\ov a \end{pmatrix},
\quad\text{where $|a|^2+|b|^2=1$}
.
$$

Denote by
 $V_n$ the space of homogeneous polynomials
of degree  $n$ in complex variables  $z_1$, $z_2$.
Define an inner product in $V_n$ assuming that
 monomials $z_1^i z_2^j$ are pairwise orthogonal and
\begin{equation}
\|z_1^i z_2^j\|^2=i!\,j!
\label{eq:zz}
\end{equation}
The group  $\SU(2)$ acts in  $V_n$ by substitutions: 
\begin{equation}
T\begin{pmatrix}a&b\\ -\ov b&\ov a \end{pmatrix} f(z_1,z_2)=
f(z_1 a-z_2 \ov b, z_1 b+z_2\ov a),
\label{eq:zamena}
\end{equation}
We write such transformations as
\begin{equation}
T(g)f(z)=f(zg)
,\end{equation}
where $z=\begin{pmatrix}z_1&z_2 \end{pmatrix}$
is a row matrix,  $zg$ is the product of the row
matrix and the square matrix.
It can be easily checked that the operators
 $T(\cdot)$ are unitary and determine a unitary representation%
\footnote{On unitary representations of  $\SU(2)$, see, for instance, 
\cite{Zhe}, \cite{Vil}; the book \cite{VK} contains numerous formulae 
with triple tensor products.}
of
$\SU(2)$.

Such representations of $\SU(2)$ are irreducible and exhaust
all irreducible representations of
  $\SU(2)$.

\SS


{\it \bf\punct Invariant vectors in tensor products.}
Consider the tensor product
$V_n\otimes V_m\otimes V_k$ of three representations
of  $\SU(2)$.
It contains an  $\SU(2)$-invariant vector iff

\SS

--- $n+m+k$ is even;

\SS

--- the numbers $n$, $m$, $k$ satisfy the 
triangle inequality,
$$
n\le m+k,\qquad m\le n+k,\qquad k\le n+m
.
$$

Recall how to construct this vector.
The space   $V_n\otimes V_m\otimes V_k$ can be regarded
as a space of polynomials in 6 variables
 $x_1$, $x_2$, $y_1$, $y_2$, $z_1$, $z_2$,
having degree  $n$ in $x_1$, $x_2$, degree 
$m$ in $y_1$, $y_2$ and degree   $k$ in $z_1$, $z_2$.
Choose the numbers $\alpha$, $\beta$, $\gamma$ from the condition
$$
\alpha+\beta=n,\qquad \alpha+\gamma=m,\qquad \beta+\gamma=k
.$$
Automatically, they are integer and nonnegative.
Then the $\SU(2)$-invariant vector 
in  $V_n\otimes V_m\otimes V_k$
is given by the formula
$$
\Xi=\Xi\bigl[\alpha,\beta,\gamma| x,y,z\bigr]=
(x_1 y_2-x_2y_1)^\alpha (y_1 z_2-y_2 z_1)^\gamma(z_1 x_2-z_2 x_1)^\beta
.
$$


{\bf\punct Statement of the problem.} Denote  
$$
G:=\SU(2)\times\SU(2)\times\SU(2)
.$$
By $K:=\SU(2)$ we denote its diagonal subgroup, i.e.,
the group of triples
 $(A,A,A)\in G$, where $A$ 
ranges in $\SU(2)$. The group  $G$ acts in 
$V_n\otimes V_m\otimes V_k$, let $\Xi$ be the vector
invariant with respect to 
 $K$.

We are going to write out (as explicit as we can)
{\it spherical functions}, i.e., the inner products
$$
\Phi_{\alpha,\beta,\gamma}(A,B,C)
:=\bigl\la T(A)\otimes T(B)\otimes T(C)\,\Xi,\Xi
\bigr\ra_{V_n\otimes V_m\otimes V_k},
\quad\text{where $(A,B,C)\in G$.} 
$$

Evidentely, the function $\Phi$ is
 $K$-biinvariant, i.e., it is constant on double cosets
$K\setminus G/K$.

By the usual orthogonality relations, the functions
$\Phi$ are pairwise orthogonal on $G$ 
with respect to the Haar measure.%
\footnote{Our space $G/K$ is a special case
of spherical varieties defined by K\"ramer \cite{Kra},
his paper initiated an extensive literature, see, for instance,
  \cite{Kno}.}

\SS


{\bf\punct Spectral curve.} To each triple  
$A$, $B$, $C\in\SU(2)$ of matrices we assign
 a {\it spectral curve} on the projective
plane  $\lambda:\mu:\nu$, this curve is defined
by the equation
$$
\det(\lambda A+\mu B+\nu C)=0
.
$$ 
Evidentely, this curve depends not on a triple
$(A,B,C)\in G$, but on а double coset
 $\in K\setminus G/K$. It can be readily checked
that a double conjugacy class can be restored from a
spectral curve.

\SS

{\it Therefore we can regard
functions $\Phi$ as functions on the space of spectral curves.}

 \SS                                                               

Since $\det A=\det B=\det C=1$,
it follows that a spectral curve has 
the form  
\begin{equation}
\lambda^2+\mu^2+\nu^2+2p\,\lambda\mu+2q\,\lambda\nu+2r\,\mu\nu=0
\label{eq:pqr}
.\end{equation}


{\bf \punct  Measure on the space of spectral curves.}

\begin{proposition}
A curve (\ref{eq:pqr}) is spectral iff  $(p,q,r)$ 
is contained in the following body 
$\Delta\subset\R^3$
\begin{eqnarray}
-1\le p\le 1,\quad -1\le q\le1,\quad -1\le r\le 1,
\\
\qquad 1-p^2-q^2-r^2+2pqr\ge 0
.\end{eqnarray}
\end{proposition}

\begin{figure}
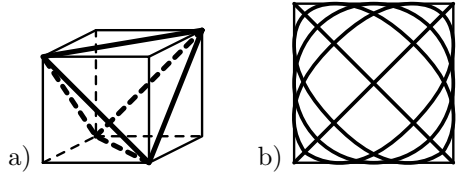


a) $\epsfbox{spectral.1}$\qquad  b) $\epsfbox{spectral.2} $

\caption{a) The body $\Delta$ is inscribed into
the unit cub $[-1,1]^3$.  
  It is tangent to the surface of the cube
along 6 segments. Namely, four vertices
$(1,1,1)$,  $(-1,-1,1)$, $(-1,1,-1)$, $(1,-1,-1)$
of the cube are conic singularities
of the surface   $1-p^2-q^2-r^2+2pqr=0$.
\newline
b) Horizontal sections of  $\Delta$ are ellipses
whose axes lies on the lines  
 $p=\pm q$.}

\end{figure}

Note that 
$$
\det\begin{pmatrix}1&p&q\\
p&1&r\\
q&r&1\end{pmatrix}=1-p^2-q^2-r^2+2pqr
.$$

Next, we have the map  $G\to\Delta$.

\begin{proposition}
The direct image of the Haar measure under the map
 $G\to\Delta$
 is the Lebesgue measure (up to a scalar factor).
\end{proposition}

By the construction, the spherical functions
 $\Phi_{\alpha,\beta,\gamma}$
constitute an orthogonal basis in 
$L^2(\Delta)$. In fact, we get a system of orthogonal
polynomials. 

\SS


{\bf\punct Generating function for spherical
functions.}

\begin{theorem}
\begin{multline*}
\sum\limits_{\alpha\ge 0,\,\beta\ge 0,\gamma \ge 0}
\frac{u^\alpha v^\beta w^\gamma}
{(\alpha!\,\beta!\,\gamma!)^2}
\Phi_{\alpha,\beta,\gamma}(p,q,r)=
\\
\Bigl(
\bigl[1+2p(vw-u)+2q(uw-v)+2r(uv-w)+u^2+v^2+w^2\bigr]^2
-uvw(1-p^2-q^2-r^2+2pqr)\Bigr)^{-1/2}
\end{multline*}
\end{theorem}

In particular,  
\begin{multline*}
\sum\limits_{\alpha\ge 0,\,\beta\ge 0,\gamma \ge 0}
\frac{u^\alpha v^\beta w^\gamma}
{(\alpha!\,\beta!\,\gamma!)^2}
\Phi_{\alpha,\beta,\gamma}(1,1,1)
=(1-u-v-w)^{-2}
=\\=
\sum
\frac{u^\alpha v^\beta w^\gamma (\alpha+\beta+\gamma+1)!}
{\alpha!\,\beta!\,\gamma!}
\end{multline*}
or
$$
\|\Xi[\alpha,\beta,\gamma;\cdot]\|^2%
_{V_\alpha\otimes V_\beta\otimes V_\gamma}
=\Phi_{\alpha,\beta,\gamma}(1,1,1) =(\alpha+\beta+\gamma+1)!
\alpha!\,\beta!\,\gamma!
$$


{\bf\punct The structure of the paper.}  Section 2
 contains preliminaries;
in \S 3 we prove the statements formulated above.
In \S 4, we get an analog of the last theorem for
an arbitrary product of groups%
\footnote{Multiple tensor products of representations
of  $\SU(2)$
were considered by numerous authors; I propose
an accidental collection of references,
\cite{How}, \cite{Kly}, \cite{Var}, \cite{Ros}.} 
$\SU(2)\times\dots\times\SU(2)$.


\section{Preliminaries. Fock space}

\COUNTERS

{\bf\punct The Fock space.}
The Fock space%
\footnote{For details, see \cite{Ner}.} 
$\F_n$ consists of entire functions on 
$\C^n$, satisfying the condition
$$
\int_{\C^n} |f(z)|^2\,e^{-|z|^2}\,dz\,d\ov z<\infty
.$$
The inner product in $\F_n$ is defined by
$$
\la f,g\ra=\int_{\C^n} f(z)\ov{g(z)}\,e^{-|z|^2}\,dz\,d\ov z
.$$


{\bf\punct The standard basis.} Monomials $z_1^{k_1}\dots z_n^{k_n}$
form an orthogonal basis $\F_n$
and
$$
\| z_1^{k_1}\dots z_n^{k_n}\|^2=k_1!\dots k_n!
$$


{\bf\punct Gauss vectors.} Let $A$ be a symmetric complex
matrix of order  $n$ satisfying $\|A\|<1$ (here a 'norm' is
the norm of linear operator in Euclidian space  $\C^n$).
Then the function
$$
\bfb[A](z)=\exp\Bigl\{\frac 12 zAz^t\Bigr\}
$$
is contained in $\F_n$. Moreover, 
\begin{equation}
\la \bfb[A],\bfb[B]\ra=\det\Bigl[(1-A\ov B)^{-1/2}\Bigr]
\label{eq:AB} 
.\end{equation}


{\bf\punct The action of $\SU(2)$ in $\F_2$.}
The group  $\SU(2)$ acts in   $\F_2$ by rotations,
$$
T\begin{pmatrix}a&b\\ -\ov b&\ov a \end{pmatrix} f(z_1,z_2)=
f(z_1 a-z_2 \ov b, z_1 b+z_2\ov a)
$$
or (in a short notation)
$$
T(g)f(z)=f(zg)
.
$$
Comparing these two formulae with (\ref{eq:zamena}), (\ref{eq:zz}),
we get
$$
\F_2=\bigoplus_{j=0}^\infty V_j
$$

\section{Calculation}

\COUNTERS

{\bf\punct  Evaluation of the radial part of the Haar measure.}
Let $(A,B,C)\in G$. Multiplying this triple
by  $A^{-1}$ from the left,
we reduce it to the form
$$(1,A^{-1}B,A^{-1}C).$$
Next, by conjugation of the second term,
we reduce the triple to the form 
\begin{equation}
\left(\begin{pmatrix} 1&0\\0&1\end{pmatrix},
\,
\begin{pmatrix} e^{i\phi}&0\\0&e^{-i\phi}\end{pmatrix},
\,
\begin{pmatrix} a&b\\ \ov b&a\end{pmatrix}\right)
\label{eq:can}
.
\end{equation}
Denote by $\wt\Delta$ the set of such triples. 
It can be readily checked that the image of the Haar measure
under the map  $G\to\wt \Delta$ is
$$
\sin^2\phi\, d\phi\times\Bigl\{\text{Haar measure onа $\SU(2)$}\Bigr\}
.$$
Next, set $b=\rho e^{i\theta}$.
In these coordinates, the Haar measure
on  $\SU(2)$ has the form 
\begin{equation}
da\,d\ov a\,d\theta,\qquad |a|\le 1,\,0<\theta<2\pi
\label{eq:haar}
.\end{equation}

The spectral curve is defined by the equation
$$
(\lambda+\mu e^{i\psi}+\nu a)
(\lambda+\mu e^{-i\psi}+\nu \ov a)+b\ov b \nu^2=0
.$$
Therefore,
$$
p=\cos\phi,\qquad q=\Re a,\qquad r=\Re(ae^{-i\phi})
.$$
We observe that the coordinate  $\theta$ does not take part
in these formulae.
The Jacobian of the pass from  $(\phi,a,\ov a)$ to
к $(p,q,r)$ is
$$
\sin^{-2}\phi=(1-p^2)^{-1}
.$$
Therefore, the image of the Haar measure is
 $dp\,dq\,dr$.

On the other hand,                                  
 $1-|a|^2\ge 0$.
Writing out this expression in new coordinates, we get
$$
1-|a|^2= 
\frac{1-p^2-q^2-r^2+2pqr}{1-p^2} 
.$$
Taking in account $-1\le p\le 1$, we get  $1-p^2-q^2-r^2+2pqr\ge 0$.

This proves Propositions 1.1 and 1.2.


\SS

{\bf\punct  Evaluation of the generating function.}
Consider the Fock space  $\F_6$, consisting of holomorphic functions
 of variables
 $x_1$, $x_2$, $y_1$, $y_2$, $z_1$, $z_2$,
$$
\F_6=\bigoplus_{n\ge 0,\, m\ge 0,\,k\ge 0}
V_n\otimes V_m\otimes V_k 
.$$
Consider the generating function
\begin{multline}
S(t,s,\sigma|x,y,z):=
\sum\limits_{\alpha\ge 0,\,\beta\ge 0,\,\gamma\ge 0}
\frac{t^\alpha s^\beta\sigma^\gamma}{\alpha!\,\beta!\,\gamma!}
\Xi[\alpha,\beta,\gamma|x,y,z]
=\\=
\exp\Bigl\{t(x_1 y_2-x_2 y_1)+s(z_1 x_2-z_2 x_1)+
\sigma(y_1z_2-y_1 z_2)  \Bigr\}
=\bfb[Q(t,s,\sigma)]
.\end{multline}
The expression in  the right-hand side  
is a Gaussian $\in\F_6$,
the matrix $Q$
has the form
$$
Q(s,t,\sigma)=
\begin{pmatrix}0&0&0&t&0&-s\\
                   0&0&-t&0&s&0\\
                   0&-t&0&0&0&\sigma\\
                   t&0 &0&0&-\sigma&0\\
                   0&s&0&-\sigma&0&0\\ 
                  -s&0&\sigma&0&0&0                    
  \end{pmatrix}
.
$$
We apply the operator
$$
T(A)\otimes T(B)\otimes T(C)f(x,y,z)=f(xA,yB,zC)
,$$
to the expression  $S(t,s,\sigma|x,y,z)\in \F_6$
and evaluate the following inner product in
 $\F_6$ 
$$
\bigl\la T(A)\otimes T(B)\otimes T(C)\, S(t,s,\sigma|x,y,z), 
S(t',s',\sigma'|x,y,z)\bigr\ra_{\F_6}
$$ 
in two ways.

First,
writing $S$ a series, we get  
$$
\la\dots,\,\dots\ra=
\sum\limits_{\alpha\ge 0,\,\beta\ge 0,\,\gamma\ge 0}
\frac{(t\ov t')^\alpha (s\ov s')^\beta(\sigma\ov\sigma')^\gamma}
{(\alpha!\,\beta!\,\gamma!)^2}\Phi_{\alpha,\beta,\gamma}(A,B,C)
.
$$

Next, denote by  $U$ the following  matrix of the order $6=2+2+2$,
$$
U=\begin{pmatrix}A&0&0\\0&B&0\\0&0&C \end{pmatrix}
.
$$
Then
\begin{multline*}
T(A)\otimes T(B)\otimes T(C)\, S(t,s,\sigma|x,y,z)=
T(A)\otimes T(B)\otimes T(C)\,\bfb[Q(s,t,\sigma)]
=\\=
\bfb[U Q(s,t,\sigma) U^t]
\end{multline*}
We can evaluate our inner product as an inner product
of Gaussian vectors by  formula 
 (\ref{eq:AB}),
$$
\la\dots,\dots\ra=
\det\Bigl(1-UQ(s,t,\sigma)U^tQ(s',t',\sigma')\Bigr)^{-1/2}
.
$$
In fact, we problem is reduced to an explicit evaluation
of a
$6\times 6$ determinant. It is sufficient to consider
 triples  (\ref{eq:can}). I can not 'refine'
this calculation, the final formula is given in Theorem 1.3;
the parameters are re-denoted as
$$
u:=s\ov s',\qquad v:=s\ov s',\qquad w:=t\ov t'
.$$


\section{Multiple tensor products}

\COUNTERS

{\bf\punct Notation.}
Now consider a product of $l$ groups,
$$
G_l=\SU(2)\times\dots\times\SU(2)
.$$
We write elements of  $G$ as
$$
(A_1,\dots, A_l),\quad\text{where $A_j\in\SU(2)$.}
.$$ 

By  $K$ we denote the diagonal subgroup, it consists of collections 
 $(A,\dots,A)$.

We consider tensor products
\begin{equation}
V_{n_1}\otimes V_{n_2}\otimes\dots \otimes V_{n_l}
\label{eq:otimes}
.\end{equation}

By $u_1$, $u_2$,\dots, $u_l$ we denote two-dimensional
complex vectors,
$$
u_j=(u_j^{(1)},u_j^{(2)})
,
$$
Thus,  $(u_1,\dots,u_l)$ is a vector $\in\C^{2l}$.
We realize the tensor product  (\ref{eq:otimes})
as the space of polynomials, whose degree  of homogeneity
in $u_j\in \C^2$
is
$n_j$ . The group $G_l$
acts in this space in the usual way,
$$
T(A_1)\otimes \dots\otimes T(A_l)f(u_1,\dots,u_l)=
f(u_1A_1,\dots,u_lA_l).
$$


{\bf\punct Invariants in tensor products.}
Let $i<j$.
Denote by $\xi_{ij}$ the following quadratic polynomial
in the variables $u$,
$$
\xi_{ij}(u):=u_i^{(1)}u_j^{(2)}- u_i^{(2)}u_j^{(1)}
.$$
Consider a symmetric  $l\times l$ matrix
$$
\alpha=\bigl\{\alpha_{ij}\bigr\}
,
$$
composed of non-negative integers, we   assume $\alpha_{ii}=0$.
By $\Xi[\alpha]$ we denote the following polynomial
$$
\Xi[\alpha](u):=\prod_{i>j} \xi_{ij}(u)^{\alpha_{ij}}
.
$$
Obviously, this polynomial is $K$-invariant. By the definition, 
$$
\Xi[\alpha]\in V_{n_1}\otimes V_{n_2}\otimes\dots \otimes V_{n_l},
\qquad\text{where $n_i=\sum_j\alpha_{ij}$}
.
$$

If matrices $\alpha$, $\beta$  satisfy the conditions 
$$
\sum_j \alpha_{ij} =n_i, \qquad \sum_j \beta_{ij} =n_i
,
$$
then we can write "spherical functions"
$$
\Phi[\alpha,\beta](A_1,\dots,A_l)
:=\bigl\la \left(\bigotimes T(A_j)\right)\Xi[\alpha],\Xi[\beta]\bigr\ra
$$

In this case a $K$-invariant vector in 
$V_{n_1}\otimes V_{n_2}\otimes\dots \otimes V_{n_l}$
is (generally speaking) nonunique; a reasonable  analog
of spherical functions is the matrix
 $\Phi[\alpha,\beta]$.

\SS


{\bf\punct Fock space.} Notice that
the Fock space  $\F_{2l}$ is a direct sum 
$$
\F_{2l}=\bigoplus_{n_1\ge 0,\dots, n_l\ge0}
 V_{n_1}\otimes V_{n_2}\otimes\dots \otimes V_{n_l}
.
$$
We can apply the same trick as in Subsection 3.2.

\SS


{\bf\punct Generating function.} 
Denote by $T_{ij}$ the following matrix of order  2 
$$
T_{ij}=\begin{pmatrix}0&t_{ij}\\-t_{ij}&0\end{pmatrix},\quad
\text{where $i<j$},
$$
$T_{ii}=0$, $T_{ji}:=-T_{ij}$. Compose a symmetric 
 $2l\times 2l$-matrix $T$ by
$$
T:=\begin{pmatrix} T_{11}&T_{12}&\dots\\
                   T_{21}&T_{22}&\dots\\
                   \vdots&\vdots&\ddots
   \end{pmatrix} 
.
$$
Next, denote by  $U$ the following block matrix  
$$
U=
\begin{pmatrix}
A_1&0&\dots\\
0&A_2&\dots\\
\vdots&\vdots&\ddots
\end{pmatrix},\quad \text{where $A_j\in \SU(2)$}
.
$$

\begin{theorem}
\begin{multline*}
\sum\limits_{n_1\ge 0,\dots,n_l\ge 0}
\,\,
\sum\limits_{\alpha:\, \sum_j \alpha_{ij}=n_i }
\,\,
\sum\limits_{\beta:\, \sum_j \beta_{ij}=n_i }
\,\,
\Bigl(
\prod_{ij} \frac{t_{ij}^{\alpha_{ij}} {\ov t'_{ij}}^{\beta_{ij}}}
             {\alpha_{ij}!\,\beta_{ij}!}\Bigr)\Phi[\alpha,\beta]
=\\=
\det(1-UTU^t T')^{-1/2}
.
\end{multline*}
\end{theorem}

In fact, the proof is contained in Subsection
3.2.

\SS

Certainly, we can remove one summation
and write the left hand-side of the identity
by 
$$
\sum\limits_{\alpha}
\,\,
\sum\limits_{\beta}
\,\,
\Bigl(
\prod_{ij} \frac{t_{ij}^{\alpha_{ij}} {\ov t'_{ij}}^{\beta_{ij}}}
             {\alpha_{ij}!\,\beta_{ij}!}\Bigr)\Phi[\alpha,\beta]
.
$$
For superfluos summands 
 $\Phi[\alpha,\beta]$ vanishes,
the corresponding powers
$t_{ij}^{\alpha_{ij}} {\ov t'_{ij}}^{\beta_{ij}}$
are absent in the series in the right-hand side.


\SS

{\bf\punct Scalar products of invariant expressions.}
Now we are going to write a generating function for
expressions
$$
\la \Xi[\alpha],\Xi[\beta] \ra=\Phi[\alpha,\beta](1,\dots,1)
$$
Obviously, it is given by the formula $\det(1-TT')^{-1/2}$.
Funnily, {\it the polynomial $\det(1-TT')$ admits a root
of forth degree}.

Denote
$$
\bfT:=\begin{pmatrix}
0&t_{12}&t_{13}&\dots
\\
-t_{12}&0&t_{23}&\dots
\\
-t_{13}&t_{23}&0&\dots
\\
\vdots&\vdots&\vdots&\ddots
\end{pmatrix}
$$

\begin{proposition}
\begin{multline*}
\sum\limits_{n_1\ge 0,\dots,n_l\ge 0}
\,\,
\sum\limits_{\alpha:\, \sum_j \alpha_{ij}=n_i }
\,\,
\sum\limits_{\beta:\, \sum_j \beta_{ij}=n_i }
\,\,
\Bigl(
\prod_{ij} \frac{t_{ij}^{\alpha_{ij}} {\ov t'_{ij}}^{\beta_{ij}}}
             {\alpha_{ij}!\,\beta_{ij}!}\Bigr)
\Phi[\alpha,\beta](1,\dots,1)
=\\=
\det(1+\bfT\bfT')^{-1}
=\det\begin{pmatrix}\bfT& 1\\-1&\bfT'\end{pmatrix}^{-1}
.
\end{multline*}
\end{proposition}

The last matrix is skew-symmetric, therefore its determinant
admits a square root (a Pfaffian).

\SS

{\sc Proof.} Consider the block diagonal matrix
  $\Sigma$ having blocks 
 $\begin{pmatrix}0&1\\1&0\end{pmatrix}$ on the diagonal. 
Then
$$
\det(1-TT')=\det\bigl(1-(T\Sigma)(\Sigma T')\bigr)
=
\det
\begin{pmatrix}  
T\Xi&1\\-1&-\Xi S
\end{pmatrix}
$$
The last determinant splits into a product of two
determinants, that are given by
$\det\begin{pmatrix}\bfT& 1\\-1&\bfT'\end{pmatrix}$.

{\tt Math.Dept., University of Vienna,

 Nordbergstrasse, 15,
Vienna, Austria

\&

Institute for Theoretical and Experimental Physics,

Bolshaya Cheremushkinskaya, 25, Moscow 117259,
Russia

\&

Mech.Math. Dept., Moscow State University,
Vorob'evy Gory, Moscow, Russia

\&

e-mail: neretin(at) mccme.ru

URL:www.mat.univie.ac.at/$\sim$neretin

wwwth.itep.ru/$\sim$neretin
}

\end{document}